\documentclass[twoside]{article}
\usepackage{amssymb}
\usepackage{indentfirst}
\usepackage{dsfont}
\usepackage{ntheorem}
\usepackage{amsmath}

\begin{document}
\baselineskip 13.5pt \noindent \thispagestyle{empty}

\markboth{\centerline{\rm Z. Tu \; \& \; S. Zhang}}{\centerline{\rm
The Schwarz Lemma at the Boundary }}

\begin{center} \large{\bf  The Schwarz Lemma at the Boundary of the Symmetrized Bidisc }
\end{center}

\begin{center}
\noindent\text{Zhenhan TU,  \; Shuo ZHANG$^*$ }\\
\vskip 4pt \noindent\small {School of Mathematics and Statistics,
Wuhan University, Wuhan, Hubei 430072, P.R. China}\\
\vskip 4pt \noindent\text{Email: zhhtu.math@whu.edu.cn (Z. Tu),\; zs.math@whu.edu.cn (S. Zhang)}
\renewcommand{\thefootnote}{{}}
\footnote{\hskip -16pt {$^{*}$Corresponding author. \\ } }
\end{center}

\begin{center}
\begin{minipage}{13cm}
{\bf Abstract.} {\small
The symmetrized bidisc $\emph{\textbf{G}}_{2}$ is defined by
$$\emph{\textbf{G}}_{2}:=\{(z_1+z_2,z_1z_2)\in\mathbb{C}^2: |z_1|<1,|z_2|<1,\; z_1,z_2\in\mathbb{C}\}.$$
It is a bounded inhomogeneous pseudoconvex domain without $\mathcal{C}^1$ boundary, and especially the symmetrized bidisc hasn't any strongly pseudoconvex boundary point and
the boundary behavior of both Carath\'eodory and Kobayashi metrics over the symmetrized bidisc is hard to describe precisely.
In this paper, we study the boundary Schwarz lemma for holomorphic self-mappings of the symmetrized bidisc $\emph{\textbf{G}}_2$, and our boundary
Schwarz lemma in the paper differs greatly from the earlier related results.
 }
\end{minipage}
\end{center}

\numberwithin{equation}{section}
\def\theequation{\arabic{section}.\arabic{equation}}

\newtheorem{lem}{Lemma}[section]
\newtheorem{thm}{Theorem}[section]
\newtheorem{prop}{Proposition}[section]
\newtheorem{cor}{Corollary}[section]

\section{Introduction}

Denote by $\emph{D}$ the unit disk in the complex plane $\mathbb{C}$. The classical Schwarz lemma
says the following:

\begin{thm}
If $f:\; \emph{D}\rightarrow \emph{D}$ is a holomorphic function that fixes the origin $0$, then
$|f(z)| \leq |z|$ for all $z\in \emph{D}$.
\end{thm}

Now the classical Schwarz lemma has become a crucial theme in many branches of mathematical research (see, for instance, Ahlfors \cite{Ah}, Rodin \cite{R}, Tsuji \cite{T} and Yau \cite{Y}).
Also, there are many results (see, for instance,  Alexander \cite{Alex}, Migliorini-Vlacci \cite{Mig},  Pommerenke \cite{Pommer}, Tauraso-Vlacci \cite{Ta}) concerning the boundary behavior of various maps, and
it is natural to consider various boundary version of the classical Schwarz lemma. There is a classical Schwarz lemma at the boundary as follows (see, for instance, Garnett \cite{Garnett}):

\begin{thm} \label{Th:1.2} Let $f:\; \emph{D}\rightarrow \emph{D}$ be a holomorphic function. If $f$ is holomorphic at $z = 1$ with $f(0)=0$ and $f(1)=1$, then $f'(1) \geqslant 1$. Moreover, the inequality
is sharp.
\end{thm}

If the assumption $f(0)=0$ in Theorem \ref{Th:1.2} is removed, then one has the following estimate $$f'(1)\geqslant\frac{|1-\overline {f(0)}|^2}{1-|f(0)|^2}>0$$ by applying Theorem \ref{Th:1.2} to the holomorphic function $g(z)=\frac{1-\overline{f(0)}}{1-f(0)} \frac{f(z)-f(0)}{1-\overline{f(0)}f(z)}$.

\vskip 5pt

The idea of Schwarz lemmas at the boundary of the unit disk has seen considerable activity in the
past 10 years or so (see, for instance,  Chelst \cite{C}, Krantz \cite{K} and Osserman \cite{O}).
Wu \cite{W} generalized the classical Schwarz lemma for holomorphic mappings
to higher dimension as follows.

\begin{thm} \label{Th:1.3}
\mbox{\rm (see \cite{W})} $($Carath\'eodory-Cartan-Kaup-Wu Theorem$)$ Let $\Omega$ be a
bounded domain in $\mathbb{C}^n$, and let $f$ be a holomorphic self-mapping of $\Omega$ which fixes a
point $p\in \Omega$. Define $J_f(p)$ as the complex Jacobian matrix of $f$ at $p$.  Then

(1) The eigenvalues of $J_f(p)$ all have modulus not exceeding 1;

(2) $|\det J_f(p)| \leqslant 1$;

(3) If $|\det J_f(p)|=1$, then $f$ is a biholomorphism of $\Omega$.
\end{thm}

Burns-Krantz \cite{BK} established a new Schwarz lemma at
the boundary, where they obtained a new rigidity result for holomorphic mappings (see Bracci-Tauraso-Vlacci \cite{BR}, Gentili-Vlacci \cite{Gentili} and Huang-Krantz \cite{HK} for related research).  Huang  \cite{H} further strengthened the result of Burns-Krantz for holomorphic mappings
with an interior fixed point. By using the boundary behavior of the Carath\'edory and Kobayashi
metrics on bounded strongly pseudoconvex domains with smooth boundary (see Graham  \cite{Graham}), recently,
Liu-Tang \cite{LT1, LT3} generalize the boundary Schwarz lemma to strongly pseudoconvex domains in $\mathbb{C}^n$.

\par  The domain $\emph{\textbf{G}}_2\subset\mathbb{C}^2$ defined by $$\emph{\textbf{G}}_2=\{(z_1+z_2,z_1z_2)\in\mathbb{C}^2 :z_1,z_2\in\emph{D}\}$$ is called {the symmetrized bidisc}. The symmetrized bidisc is a bounded pseudoconvex domain. It is important because it is the first known example of a bounded pseudoconvex
domain for which the Carath\'eodory and Lempert functions coincide, but which cannot be exhausted by
domains biholomorphic to convex ones (see Costara \cite{Co}).  Moreover, the symmetrized bidisc plays also an important role in solving the
Pick-Nevanlinna Interpolation Problem in dimension two (cf.  Agler-Young \cite{AY}). The symmetrized bidisc has been recently studied by many authors, e.g., Agler-Lykova-Young \cite{ALY}, Agler-Young \cite{You}, Frosini-Vlacci \cite{FV2}, Jarnicki-Pflug \cite{JP},  Su-Tu-Wang \cite{S}  and Trybu{\l}a \cite{MT}.
Specially, Agler-Young \cite{You2} gave a Schwarz lemma for the symmetrized bidisc in 2001.

Following this line, we study the boundary Schwarz lemma for the symmetrized bidisc $\emph{\textbf{G}}_2$ in this paper.
Note that the symmetrized bidisc is a bounded inhomogeneous domain without smooth boundary, and especially the symmetrized bidisc has no strongly pseudoconvex boundary point
and the boundary behavior of both Carath\'eodory and Kobayashi metrics over the symmetrized bidisc is hard to describe precisely. We need to find a different approach for such
a study. Because the symmetrized bidisc has no strong pseudoconvex boundary point, our boundary
Schwarz lemma in the paper differs greatly from the earlier related results (e.g., see Liu-Tang \cite{LT1, LT3}).

\section{Preliminaries}

In this section,we exhibit some notations and collect several basic lemmas, which will be used in the subsequent section.

\begin{lem}$(\mathrm{See}\,\, \emph{\cite{You}})$ \label{Lm:2.1}
Let $s, p\in \mathbb{C}$. The following statements are equivalent:\par $(1)$ $(s, p)\in\emph{\textbf{G}}_2 $; \par $(2)$ the roots of the equation $z^2-sz+p=0$ lie in $D$; \par $(3)$ $|s-\overline{s}p|<1-|p|^2$; \par $(4)$ $|s|<2$ and for all $\omega\in \mathbb{T}$,$$\Big|\frac{2p-\omega s}{2-\overline{\omega}s}\Big|<1;$$ \par $(5)$ $|p|<1$ and there exists $\beta\in D$,such that $s=\beta p +\overline{\beta}$; \par $(6)$ $2|s-\overline{s}p|+|s^2-4p|+|s|^2<4$.
\end{lem}

\par Suppose that $\Omega$ is a domain in $\mathbb{C}^n$. Let $H(\Omega ,\emph{D})$ be the set of all holomorphic mappings from $\Omega$ into $\emph{D}$. For any $z\in\Omega$ and $\xi\in\mathbb{C}^n$,$$F_C ^\Omega(z, \xi)=\sup\{|\emph{J}_f(z)\xi|: f\in H(\Omega ,\emph{D}), f(z)=0\}$$ is said to be the infinitesimal form of Carath\'eodory metric of $\Omega$. Here $\emph{J}_f(z)=(\frac{\partial f}{\partial z_1},\cdots, \frac{\partial f}{\partial z_n})$.

Agler-Young \cite{You} give the form of the Carath\'eodory metric of the symmetrized bidisc $\emph{\textbf{G}}_2$.

\begin{lem}$(\mathrm{See}\,\, \emph{\cite{You}})$ \label{Lm:2.2}
If $z=(s, p)\in \emph{\textbf{G}}_2 $ and $\xi=(\xi_1, \xi_2)'\in \mathbb{C}^2$, let $\mathbb{T}$ be the unit circle, then
$$F_C(z, \xi)=\sup_{\omega\in\mathbb{T}}\bigg|\frac{\xi_1(1-\omega^{2}p)-\xi_2(2-\omega s)\omega}{(s-\overline{s}p)\omega^2-2(1-|p|^2)\omega+\overline{s}-\overline{p}s}\bigg|.$$
\end{lem}
\par From the proof of Corollary 4.4 in \cite{You}, we can get another formula of $F_C(z, \xi)$ as $$F_C(z, \xi)=\sup_{\omega\in\mathbb{T}}\frac{2|\xi_1(1-\omega^{2}p)-\xi_2(2-\omega s)\omega|}{|2-\omega s|^2-|2\omega p-s|^2},$$
which is very useful for our calculation.

The royal variety $\varSigma$ of the symmetrized bidisc $\emph{\textbf{G}}_2$ plays an important role in the study of the symmetrized bidisc. The royal variety of $\emph{\textbf{G}}_2$ is defined by
$$\varSigma:=\{(2\lambda,\lambda^2): \lambda\in D\}\subset\emph{\textbf{G}}_2.$$
Note that Jarnicki-Pflug \cite{JP} completely describe the group of holomorphic automorphisms for the symmetrized bidisc $\emph{\textbf{G}}_2$ as follows.

\begin{lem}$(\mathrm{See}\,\, \emph{\cite{JP}})$ \label{Lm:2.3}
$$ \mathrm {Aut}(\emph{\textbf{G}}_2)=\{H_h: h\in \mathrm {Aut(D)} \},$$
where $H_h(\lambda_1+\lambda_2,\lambda_1\lambda_2)=(h(\lambda_1)+h(\lambda_2),h(\lambda_1)h(\lambda_2))$, $(\lambda_1+\lambda_2,\lambda_1\lambda_2)\in\emph{\textbf{G}}_2$ with $\lambda_1,\lambda_2\in D$. Because $\mathrm{Aut}(\emph{\textbf{G}}_2)$ does not act transitively on $\emph{\textbf{G}}_2$, the symmetrized bidisc is inhomogeneous. But the group $\mathrm{Aut}(\emph{\textbf{G}}_2)$ acts transitively on $\varSigma$.
\end{lem}

The following lemma characterizes the contraction property of the Carath\'eodory metric, which is also a version of Schwarz lemma.

\begin{lem} $(\mathrm{See}\,\, \emph{\cite{FV}})$ \label{Lm:2.4}
Let $\phi :\emph{\textbf{G}}_2 \rightarrow \emph{\textbf{G}}_2$ be a holomorphic mapping. Then for any $z\in \emph{\textbf{G}}_2$ and $\xi\in\mathbb{C}^2$,$$F_C(\phi(z),J_\phi (z)\xi)\leqslant F_C(z,\xi).$$ Moreover, if $\phi :\emph{\textbf{G}}_2 \rightarrow \emph{\textbf{G}}_2$ be a biholomorphic mapping, then we have $$F_C(\phi(z),J_\phi (z)\xi)=F_C(z,\xi).$$
\end{lem}

With the form of the Carath\'eodory metric of $\emph{\textbf{G}}_2$ in Lemma \ref{Lm:2.2} and the holomorphic automorphism group of $\emph{\textbf{G}}_2$ in Lemma \ref{Lm:2.3}, we can get the explicit formula of the Carath\'eodory metric at some points of $\emph{\textbf{G}}_2$ as follows.

\begin{lem} \label{Lm:2.5}
Let $z=(s_0, p_0)=(2\alpha, \alpha^2)\in \emph{\textbf{G}}_2$ with $\alpha\in\mathbb{C}$, $|\alpha|<1$ and $\xi=(\xi_1, \xi_2)'\in\mathbb{C}^2$, then
$$F_C(z, \xi)=\frac{|(1+|\alpha|^2)\xi_1-2\overline{\alpha}\xi_2|+2|\alpha\xi_1-\xi_2|}{2(1-|\alpha|^2)^2}.$$
\end{lem}

\noindent
{\bf Proof}.
By Lemma \ref{Lm:2.3}, we can take $H_{h_\alpha}\in \mathrm{Aut}(\emph{\textbf{G}}_2)$ such that $H_{h_\alpha}(s_0, p_0)=(0,0)$, where
$$h_\alpha(\lambda)=\frac{\lambda-\alpha}{1-\overline{\alpha}\lambda}\in\mathrm{Aut}(D) ,$$
\begin{align*}
H_{h_\alpha}(s, p)&=(h_\alpha(\lambda_1)+h_\alpha(\lambda_2), h_\alpha(\lambda_1)h_\alpha(\lambda_2))\\&=\left(\frac{\lambda_1-\alpha}{1-\overline{\alpha}\lambda_1}+ \frac{\lambda_2-\alpha}{1-\overline{\alpha}\lambda_2}, \frac{\lambda_1-\alpha}{1-\overline{\alpha}\lambda_1}\frac{\lambda_2-\alpha}{1-\overline{\alpha}\lambda_2}\right)\\&=\left(\frac{(1+|\alpha|^2)s-2\overline{\alpha}p-2\alpha}{1-\overline{\alpha}s+{\overline{\alpha}}^2p}, \frac{p-\alpha s+\alpha^2}{1-\overline{\alpha}s+{\overline \alpha}^2 p}\right)
\end{align*}
with $(s,p)=(\lambda_1+\lambda_2,\lambda_1\lambda_2)\in \emph{\textbf{G}}_2$, $\lambda_1$ and $\lambda_2\in D$.
\par A direct calculation shows that
\begin{align*}
\emph{J}_{H_{h_\alpha}}(s_0,p_0) = \frac{1}{(1-|\alpha|^2)^2}
\left( \begin{array}{ccc}
1+|\alpha|^2 & -2\overline \alpha  \\
-\alpha & 1  \\
 \end{array} \right)  .
\end{align*}
\par By Lemma \ref{Lm:2.2}, we have $$F_C((0, 0),(v_1,v_2)')=\frac{|v_1|+2|v_2|}{2}$$ for any $(v_1,v_2)'\in \mathbb{C}^2$. By Lemma \ref{Lm:2.4}, it follows that
\begin{align*}
F_C(z, \xi) & =  F_C((s_0,p_0),(\xi_1,\xi_2)')\\ & =  F_C(H_{h_\alpha}(s_0,p_0),\emph{J}_{H_{h_\alpha}}(s_0,p_0)(\xi_1,\xi_2)')\\&=F_C\left((0,0),\frac{1}{(1-|\alpha|^2)^2}
\left( \begin{array}{ccc}
1+|\alpha|^2 & -2\overline \alpha  \\
-\alpha & 1  \\
 \end{array} \right)
 \left( \begin{array}{ccc}
\xi_1 \\\xi_2  \\
 \end{array} \right)\right)\\&= F_C\left((0,0),\frac{1}{(1-|\alpha|^2)^2}((1+|\alpha|^2)\xi_1-2\overline \alpha \xi_2, -\alpha\xi_1+\xi_2)\right)\\&=\frac{|(1+|\alpha|^2)\xi_1-2\overline\alpha\xi_2|+2|\alpha\xi_1-\xi_2|}{2(1-|\alpha|^2)^2}.
\end{align*}
The proof is complete.

\begin{lem} \label{Lm:2.6}
If $s\in \mathbb{C}$ with $|s|<1$, then, for any $\beta\in \mathbb{C}$, we have $$F_C((s,0),(\beta, s\beta)')=\frac{|s|+1}{2}|\beta|.$$
\end{lem}

\noindent
{\bf Proof}.
By Lemma \ref{Lm:2.2}, for any $s\in \mathbb{C}$, $|s|<1$ and $\beta\in \mathbb{C}$, we have
\begin{align}
F_C((s,0),(\beta, s\beta)')\nonumber &=\sup_{\omega\in \mathbb{T}}\bigg|\frac{\beta-s\beta(2-\omega s)\omega}{s\omega^2-2\omega+\overline{s}}\bigg|\\ \nonumber&=|\beta|\sup_{\omega\in \mathbb{T}}\bigg|\frac{1-(2-\omega s)\omega s}{s\omega^2-2\omega+\overline{s}}\bigg|\\ \nonumber&=|\beta|\sup_{\omega\in \mathbb{T}}\bigg|\frac{s^2\omega^2-2\omega s+1}{s\omega ^2-2\omega+\overline{s}}\bigg|\\&\nonumber=|\beta|\sup_{\omega\in \mathbb{T}}\bigg|\frac{s^2\omega^2-2\omega s+|s|^2+1-|s|^2}{s\omega ^2-2\omega+\overline{s}}\bigg|\\ \label{eq2.1} &=|\beta|\sup_{\omega\in \mathbb{T}}\bigg|s+\frac{1-|s|^2}{s\omega^2-2\omega+\overline{s}}\bigg|.
\end{align}
 Notice that
\begin{eqnarray*}
   & &\bigg|s+\frac{1-|s|^2}{s\omega^2-2\omega+\overline{s}}\bigg|^2\\
   &=& \left(s+\frac{1-|s|^2}{s\omega^2-2\omega+\overline{s}}\right)\left(\overline s+\frac{1-|s|^2}{\bar s\bar{\omega}^2-2\overline{\omega}+s}\right)\\
   &=& |s|^2+(1-|s|^2)\frac{s}{\bar s\bar{\omega}^2-2\overline{\omega}+s}+(1-|s|^2)\frac{\overline{s}}{s\omega^2-2\omega+\overline{s}}+ \frac{(1-|s|^2)^2}{(s\omega^2-2\omega+\overline{s})(\bar s\bar{\omega}^2-2\overline{\omega}+s)}\\
   &=& |s|^2+(1-|s|^2)\frac{s^2\omega^2+\bar{s}^2\bar{\omega}^2-2(s\omega+\bar{s}\bar{\omega})+|s|^2+1}{s^2\omega^2+\bar{s}^2\bar{\omega}^2-4(s\omega+\bar{s}\bar{\omega})+2|s|^2+4}\\
   &=& |s|^2+\frac{1-|s|^2}{2}\left(1+\frac{s^2\omega^2+\bar{s}^2\bar{\omega}^2-2}{s^2\omega^2+\bar{s}^2\bar{\omega}^2-4(s\omega+\bar{s}\bar{\omega})+2|s|^2+4}\right).
\end{eqnarray*}
\par This, together with  \eqref{eq2.1}, implies
\begin{eqnarray}
\nonumber  & & F_C((s,0),(\beta, s\beta)')\nonumber\\
\nonumber  &=& |\beta|\sup_{\omega\in \mathbb{T}}\bigg|s+\frac{1-|s|^2}{s\omega^2-2\omega+\overline{s}}\bigg|\\ \label{eq2.2}
&=& |\beta|\sqrt{|s|^2+\frac{1-|s|^2}{2}\left[ 1+\sup_{\omega\in\mathbb{T}}\left(\frac{s^2\omega^2+\bar{s}^2\bar{\omega}^2-2}{s^2\omega^2+\bar{s}^2\bar{\omega}^2-4(s\omega+\bar{s}\bar{\omega})+2|s|^2+4}\right) \right]}.
\end{eqnarray}
So we only need to find $$M:=\sup_{\omega\in\mathbb{T}}\left(\frac{s^2\omega^2+\bar{s}^2\bar{\omega}^2-2}{s^2\omega^2+\bar{s}^2\bar{\omega}^2-4(s\omega+\bar{s}\bar{\omega})+2|s|^2+4}\right).$$
\par Let $c=s\omega$ with $c=a+bi$, then $|c|=|s|<1$ and $a^2+b^2=|s|^2$, where $a,b\in \mathbb{R}$. Notice that
\begin{align*}
\frac{s^2\omega^2+\bar{s}^2\bar{\omega}^2-2}{s^2\omega^2+\bar{s}^2\bar{\omega}^2-4(s\omega+\bar{s}\bar{\omega})+2|s|^2+4}&=\frac{c^2+\overline{c}^2-2}{c^2+\overline{c}^2-4(c+\overline{c})+2|s|^2+4}
\\&=\frac{2(2a^2-|s|^2)-2}{2(2a^2-|s|^2)-8a+2|s|^2+4}\\&=\frac{2a^2-|s|^2-1}{2(1-a)^2},
\end{align*}
where $a\in\mathbb{R}$ and $|a|\leqslant|s|$.
\par Let $f(a)=\frac{2a^2-|s|^2-1}{(1-a)^2}$, $a\in[-|s|,|s|]$.
Then a simple calculation shows that $$f'(a)<0,\,\, \,a\in[-|s|,|s|],$$
which implies $f(a)$ is decreasing on $a\in[-|s|,|s|]$. So we have
\begin{align*}
M&=\sup_{\omega\in\mathbb{T}}\left(\frac{s^2\omega^2+\bar{s}^2\bar{\omega}^2-2}{s^2\omega^2+\bar{s}^2\bar{\omega}^2-4(s\omega+\bar{s}\bar{\omega})+2|s|^2+4}\right)
\\&=\frac{1}{2}\max_{a\in[-|s|,|s|]}f(a)=\frac{1}{2}f(-|s|)=\frac{1}{2}\frac{|s|-1}{|s|+1}.
\end{align*}
This, together with \eqref{eq2.2}, implies $$F_C((s,0),(\beta, s\beta)')=\frac{|s|+1}{2}|\beta|.$$
The proof is complete.

\section{Main Results}

In this section, we present the main results in the article. We study the Schwarz lemma at the boundary points $(e^{i\theta},0)$, $(0,e^{i\theta})$ and $(2\alpha,\alpha^2)$ ($\theta\in[0,2\pi)$, $\alpha\in\mathbb{C}$ and $|\alpha|=1$), which will be exhibited in Theorem \ref{Th:3.1}, Theorem \ref{Th:3.2} and Theorem \ref{Th:3.3} respectively.
\vbox{}

\begin{thm} \label{Th:3.1}
Let $f:\emph{\textbf{G}}_2\rightarrow\emph{\textbf{G}}_2$ be a holomorphic mapping and set $f=(f_1, f_2)$. Let $z_0=(e^{i\theta},0)\in\partial\emph{\textbf{G}}_2$, where $\theta\in[0,2\pi)$. If $f$ is holomorphic at $z_0$ and $f(z_0)=z_0$, then the following statements hold:
\par $\mathrm{(\romannumeral1)}$  $\lambda:=\frac{\partial f_1}{\partial s}(z_0)-e^{-i\theta}\frac{\partial f_2}{\partial s}(z_0)$ and
$\mu:=\frac{\partial f_1}{\partial s}(z_0)+e^{i\theta}\frac{\partial f_1}{\partial p}(z_0)$ are all eigenvalues of $J_f(z_0)$;
\par $\mathrm{(\romannumeral2)}$  $\lambda\geqslant\frac{1}{2}\frac{|1-\overline{g(0)}|^2}{1-|g(0)|^2}>0$, where $g(0)=e^{-2i\theta}\frac{e^{i\theta}f_1(0, 0)-2f_2(0,0)}{2-e^{-i\theta}f_1(0,0)}$;
\par $\mathrm{(\romannumeral3)}$  The normal vector $(e^{i\theta}, -e^{2i\theta})'$ to $\partial\emph{\textbf{G}}_2$ at $z_0$ is an eigenvector of  $\overline{J_f(z)'}$ with respect to $\lambda$.
That is, $$\overline{J_f(z_0)'}
\left( \begin{array}{ccc}
e^{i\theta} \\ -e^{2i\theta}  \\
 \end{array} \right)=\lambda\left( \begin{array}{ccc}
e^{i\theta} \\ -e^{2i\theta}  \\
 \end{array} \right).$$
 Moreover, we have $$\lambda:=\frac{\partial f_1}{\partial s}(z_0)-e^{-i\theta}\frac{\partial f_2}{\partial s}(z_0)=\frac{\partial f_2}{\partial p}(z_0)-e^{i\theta}\frac{\partial f_1}{\partial p}(z_0);$$
\par $\mathrm{(\romannumeral4)}$  $|\mu|\leqslant1$. The eigenvectors of $J_f(z_0)$ with respect to $\mu$
have the form of
$$(\alpha,e^{i\theta}\alpha)' \in T_{z_0}^{(1,0)}(\partial\emph{\textbf{G}}_2 ),$$
where $\alpha\in\mathbb{C}-\{0\}$. Moreover, we have
$$\mu:=\frac{\partial f_1}{\partial s}(z_0)+e^{i\theta}\frac{\partial f_1}{\partial p}(z_0)=e^{-i\theta}\frac{\partial f_2}{\partial s}(z_0)+\frac{\partial f_2}{\partial p}(z_0);$$
\par $\mathrm{(\romannumeral5)}$ $|\det J_f(z_0)|\leqslant\lambda$ and $|\mathrm{tr}\,J_f(z_0)|\leqslant\lambda+1$.
\par Moreover, the inequalities in $\mathrm{(\romannumeral4)}$ and $\mathrm{(\romannumeral5)}$ are sharp
\end{thm}

\noindent
{\bf Proof}.
The proof is divided into four steps.\\
\emph{Step} 1.   From the equivalence of (1) and (6) in Lemma \ref{Lm:2.1}, we can see that
$$h(s,p)=2|s-\overline{s}p|+|s^2-4p|+|s|^2-4$$
is a defining function for $\emph{\textbf{G}}_2$. A simple calculation shows that
$$\frac{\partial h}{\partial s}=\frac{\overline{s}-2s\overline{p}+\overline{s}|p|^2}{|s-\overline{s}p|}+\frac{|s|^2\overline{s}-4s\overline{p}}{|s^2-4p|}+\overline{s},$$
$$\frac{\partial h}{\partial p}=\frac{-\overline{s}^2+|s|^2\overline{p}}{|s-\overline{s}p|}-\frac{2(\overline{s}^2-4\overline{p})}{|s^2-4p|}.$$
We can easily verify that $\partial\emph{\textbf{G}}_2$ has $C^1$ boundary near $z_0$ and $\frac{\partial h}{\partial s}(z_0)=3e^{-i\theta}$, $\frac{\partial h}{\partial p}(z_0)=-3e^{-2i\theta}$.
\par Denote by $T_{z_0}(\partial\emph{\textbf{G}}_2 )$ and $T_{z_0}^{(1,0)}(\partial\emph{\textbf{G}}_2 )$ the tangent space and holomorphic tangent space to $\partial\emph{\textbf{G}}_2 $  at $z_0=(e^{i\theta},0)\,(\theta\in[0,2\pi))$ respectively. Then
\begin{align*}
T_{z_0}(\partial\emph{\textbf{G}}_2 )&=\{\alpha=(\alpha_1,\alpha_2)'\in\mathbb{C}^2:\,\Re\,(\frac{\partial h}{\partial s}(z_0)\alpha_1+\frac{\partial h}{\partial p}(z_0)\alpha_2)=0\}\\&=\{\alpha=(\alpha_1,\alpha_2)'\in\mathbb{C}^2:\,\, \Re\,(e^{-i\theta}\alpha_1-e^{-2i\theta}\alpha_2)=0\},
\end{align*}
\begin{align*}
T_{z_0}^{(1,0)}(\partial\emph{\textbf{G}}_2 )&=\{\alpha=(\alpha_1,\alpha_2)'\in\mathbb{C}^2:\,\frac{\partial h}{\partial s}(z_0)\alpha_1+\frac{\partial h}{\partial p}(z_0)\alpha_2=0\}\\&=\{\alpha=(\alpha_1,\alpha_2)'\in\mathbb{C}^2:\,e^{i\theta}\alpha_1=\alpha_2\}.
\end{align*}
\par Since $f$ is holomorphic at $z_0$, we assume that $f$ is holomorphic on a neighborhood $V$ of $z_0$. For any $\alpha=(\alpha_1,\alpha_2)'\in T_{z_0}(\partial\emph{\textbf{G}}_2 )$, take the smooth curve $\gamma:[-1,1]\rightarrow\partial\emph{\textbf{G}}_2 $ such that $\gamma(0)=z_0$, $\gamma'(0)=\alpha$ and $\gamma([-1,1])\subset V$.
 Obviously, $h[f(\gamma(t))]\in C^1([-1,1])$ and $h[f(\gamma(t))]\leqslant0$, $t\in[-1,1]$. Hence
 $$\max_{t\in[-1,1]}h[f(\gamma(t))]=h[f(\gamma(0))]=h(z_0)=0.$$
 This implies that
\begin{align}
0\nonumber&=\frac{d}{dt}\left(h[f(\gamma(t))]\right)\Big|_{t=0}\\\nonumber&= \frac{\partial h}{\partial s}(z_0)\frac{df_1(\gamma(t))}{dt}\Big|_{t=0}+\frac{\partial h}{\partial p}(z_0)\frac{df_2(\gamma(t))}{dt}\Big|_{t=0}+\frac{\partial h}{\partial \overline s}(z_0)\frac{d\overline {f_1(\gamma(t))}}{dt}\Big|_{t=0}+\frac{\partial h}{\partial \overline p}(z_0)\frac{d\overline {f_2(\gamma(t))}}{dt}\Big|_{t=0}\\\nonumber&=2\Re \, (\frac{\partial h}{\partial s}(z_0)\frac{df_1(\gamma(t))}{dt}\Big|_{t=0}+\frac{\partial h}{\partial p}(z_0))\frac{df_2(\gamma(t))}{dt}\Big|_{t=0})
\\\nonumber&=2\Re\,[3e^{-i\theta}(\frac{\partial f_1}{\partial s}(z_0)\alpha_1+\frac{\partial f_1}{\partial p}(z_0)\alpha_2)-3e^{-2i\theta}(\frac{\partial f_2}{\partial s}(z_0)\alpha_1+\frac{\partial f_2}{\partial p}(z_0)\alpha_2)]
\\ \label{eq3.1}&=6\Re\,(e^{-i\theta},-e^{-2i\theta})J_f(z_0)\alpha.
\end{align}
It follows that $$J_f(z_0)T_{z_0}(\partial\emph{\textbf{G}}_2 )\subset T_{z_0}(\partial\emph{\textbf{G}}_2 ).$$
Consequently, there exists $\lambda\in \mathbb{R}$ such that $$(e^{-i\theta},-e^{-2i\theta})J_f(z_0)=\lambda(e^{-i\theta},-e^{-2i\theta}).$$
That is $$\overline{J_f(z_0)'}\left( \begin{array}{ccc}
e^{i\theta} \\ -e^{2i\theta}  \\
 \end{array} \right) =\lambda\left( \begin{array}{ccc}
e^{i\theta} \\ -e^{2i\theta}  \\
 \end{array} \right).$$
It follows that $$\lambda=\frac{\partial f_1}{\partial s}(z_0)-e^{-i\theta}\frac{\partial f_2}{\partial s}(z_0)=\frac{\partial f_2}{\partial p}(z_0)-e^{i\theta}\frac{\partial f_1}{\partial p}(z_0).$$
This means that $\lambda$ is an eigenvalue of $\overline{J_f(z_0)'}$. Since $\lambda\in\mathbb{R}$, it is also an eigenvalue of $J_f(z_0)$.
\par Let $\nabla h(z_0)=2(\frac{\partial h}{\partial \overline{s}}(z_0), \frac{\partial h}{\partial \overline{p}}(z_0))'$ be the outward normal vector to $\partial\emph{\textbf{G}}_2$ at $z_0$. Then we have $\nabla h(z_0)=(6e^{i\theta}, -6e^{2i\theta})'$. It follows that $(e^{i\theta},-e^{2i\theta})'$ is  a normal vector to $\partial\emph{\textbf{G}}_2$ at $z_0$.
The proof of $\mathrm{(\romannumeral3)}$ is complete.\\

\emph{Step} 2.   Let $$g(\xi)=e^{-2i\theta}\frac{e^{i\theta}f_1(\xi e^{i\theta},0)-2f_2(\xi e^{i\theta},0)}{2-e^{-i\theta}f_1(\xi e^{i\theta},0)},\,\,\,\xi\in D.$$
From the equivalence of (1) and (4) in Lemma \ref{Lm:2.1}, we know that $g$ is a holomorphic mapping from $D$ to $D$ and $g$ is holomorphic at $\xi=1$ with $g(1)=e^{-2i\theta}\frac{e^{i\theta}f_1(e^{i\theta},0)-2f_2(e^{i\theta},0)}{2-e^{-i\theta}f_1(e^{i\theta},0)}=1$. By
 Theorem \ref{Th:1.2}, we have
$$g'(1)=2\bigg(\frac{\partial f_1}{\partial s}(z_0)-e^{-i\theta}\frac{\partial f_2}{\partial s}(z_0)\bigg)=2\lambda\geqslant\frac{|1-\overline{g(0)}|^2}{1-|g(0)|^2}>0.$$
That is $$\lambda\geqslant\frac{1}{2}\frac{|1-\overline{g(0)}|^2}{1-|g(0)|^2}>0,$$
where $g(0)=e^{-2i\theta}\frac{e^{i\theta}f_1(0, 0)-2f_2(0,0)}{2-e^{-i\theta}f_1(0,0)}$.
The proof of $\mathrm{(\romannumeral2)}$ is
complete.
\\

\emph{Step} 3.  Notice that for any $\alpha\in T_{z_0}^{(1,0)}(\partial\emph{\textbf{G}}_2 )\subset T_{z_0}(\partial\emph{\textbf{G}}_2 )$, we have $e^{i\theta}\alpha\in T_{z_0}^{(1,0)}(\partial\emph{\textbf{G}}_2 )$. Similar to the proof of \eqref{eq3.1}, we obtain $$ \Re \,\,[(e^{-i\theta},-e^{-2i\theta})J_f(z_0)(e^{i\theta}\alpha)]=0$$ for any $\theta\in\mathbb{R}$.
Let $\omega=(e^{-i\theta},-e^{-2i\theta})J_f(z_0)\alpha$, then $\Re\,\,(e^{i\theta}\omega)=0$ for any $\theta\in\mathbb{R}$. Take $\theta=0$ and $\theta=\frac{\pi}{2}$, we obtain $\Re\,\,\omega=0$ and $\Im \,\,\omega=0$ respectively. That is $\omega=(e^{-i\theta},-e^{-2i\theta})J_f(z_0)\alpha=0$, which means $$J_f(z_0)T_{z_0}^{(1,0)}(\partial\emph{\textbf{G}}_2 )\subset T_{z_0}^{(1,0)}(\partial\emph{\textbf{G}}_2 ).$$
Hence, $ J_f(z_0)$ is a linear transformation on 1-dimensional complex vector space
$T_{z_0}^{(1,0)}(\partial\emph{\textbf{G}}_2 )=\{\alpha=(\alpha_1,\alpha_2)'\in\mathbb{C}^2:\,e^{i\theta}\alpha_1=\alpha_2\}$.
Let $\mu\in\mathbb{C}$ be the eigenvalue of the linear transformation $J_f(z_0)$ on $T_{z_0}^{(1,0)}(\partial\emph{\textbf{G}}_2 )$
 and $(\alpha,e^{i\theta}\alpha)'\in T_{z_0}^{(1,0)}(\partial\emph{\textbf{G}}_2 )$ be an eigenvector of $J_f(z_0)$ with respect
 to $\mu$, where $\alpha\in\mathbb{C}-\{0\}$. That is
 $$J_f(z_0)
\left( \begin{array}{ccc}
\alpha \\ e^{i\theta}\alpha  \\
 \end{array} \right)=\mu\left( \begin{array}{ccc}
\alpha \\ e^{i\theta}\alpha  \\
 \end{array} \right).$$
 It follows that $$\mu=\frac{\partial f_1}{\partial s}(z_0)+e^{i\theta}\frac{\partial f_1}{\partial p}(z_0)=e^{-i\theta}\frac{\partial f_2}{\partial s}(z_0)+\frac{\partial f_2}{\partial p}(z_0).$$

\par Let $t\in(0,1)$. By Lemma \ref{Lm:2.2}, we obtain
\begin{align*}
&\,\,\,\,\,\,\,\,\,F_C[f(te^{i\theta},0),J_f(te^{i\theta},0)(\alpha,te^{i\theta}\alpha)']\\
&=\sup_{\omega\in\mathbb{T}}\frac{2\bigg|J_{f_1}(te^{i\theta},0)
\left( \begin{array}{ccc}
\alpha\\ te^{i\theta}\alpha   \\
 \end{array} \right)
 (1-\omega^2 f_2(te^{i\theta},0))-J_{f_2}(te^{i\theta},0)
 \left( \begin{array}{ccc}
\alpha\\ te^{i\theta}\alpha  \\
 \end{array} \right)(2-\omega f_1(te^{i\theta},0))\omega\bigg|}{|2-\omega f_1(te^{i\theta},0)|^2-|2\omega f_2(te^{i\theta},0)-f_1(te^{i\theta},0)|^2}\\
&\geqslant  2|\alpha| \frac{\bigg|J_{f_1}(te^{i\theta},0)
\left( \begin{array}{ccc}
1\\ te^{i\theta}   \\
 \end{array} \right)(1-e^{-2i\theta}f_2(te^{i\theta},0))+J_{f_2}(te^{i\theta},0)
\left( \begin{array}{ccc}
1\\ te^{i\theta}  \\
 \end{array} \right)(2+e^{-i\theta}f_1(te^{i\theta},0))e^{-i\theta}\bigg|}{|2+e^{-i\theta} f_1(te^{i\theta},0)|^2-|2e^{-i\theta} f_2(te^{i\theta},0)+f_1(te^{i\theta},0)|^2}.
\end{align*}
By Lemma \ref{Lm:2.4}, the contraction property of the Carath\'{e}odory metric, we have
$$F_C[f(te^{i\theta},0),J_f(te^{i\theta},0)(\alpha,te^{i\theta}\alpha)']\leqslant F_C((te^{i\theta},0),(\alpha,te^{i\theta}\alpha)').$$
This, together with Lemma \ref{Lm:2.6}, yields
\begin{align*}
&\;\;\;\;\; |\alpha|\frac{t+1}{2}\\
&=F_C((te^{i\theta},0),(\alpha,te^{i\theta}\alpha)')\\
&\geqslant F_C[f(te^{i\theta},0),J_f(te^{i\theta},0)(\alpha,te^{i\theta}\alpha)']\\
&\geqslant 2|\alpha| \frac{\bigg|J_{f_1}(te^{i\theta},0)
\left( \begin{array}{ccc}
1\\ te^{i\theta}   \\
 \end{array} \right)(1-e^{-2i\theta}f_2(te^{i\theta},0))+J_{f_2}(te^{i\theta},0)
\left( \begin{array}{ccc}
1\\ te^{i\theta}  \\
 \end{array} \right)(2+e^{-i\theta}f_1(te^{i\theta},0))e^{-i\theta}\bigg|}{|2+e^{-i\theta} f_1(te^{i\theta},0)|^2-|2e^{-i\theta} f_2(te^{i\theta},0)+f_1(te^{i\theta},0)|^2}.
\end{align*}
for any $t\in(0,1)$.
\par As $t\rightarrow1^{-}$, we obtain $|\mu|\leqslant1$. Moreover, by taking the identity mapping on $\emph{\textbf{G}}_2$, it is easy to check that the inequality is sharp.
The proof of $\mathrm{(\romannumeral4)}$ is complete.
\\

\emph{Step} 4.  Now we claim that $\lambda$ and $\mu$ are the all eigenvalues of the linear transformation $J_f(z_0)$ on $\mathbb{C}^2$.
\par If $\lambda\neq\mu$, then $\lambda$ and $\mu$ are all the eigenvalues of the linear transformation $J_f(z_0)$ on $\mathbb{C}^2$.
\par If $\lambda=\mu$, consider the characteristic polynomial of  $J_f(z_0)$:
$$\det\left(xI_2-J_f(z_0)\right)=x^2-\left(\frac{\partial f_1}{\partial s}(z_0)+\frac{\partial f_2}{\partial p}(z_0)\right)x+\frac{\partial f_1}{\partial s}(z_0)\frac{\partial f_2}{\partial p}(z_0)-\frac{\partial f_1}{\partial p}(z_0)\frac{\partial f_2}{\partial s}(z_0).$$
By $\mathrm{(\romannumeral3)}$ and $\mathrm{(\romannumeral4)}$ which we have proved before, together with $\lambda=\mu$, we obtain
$$\frac{\partial f_1}{\partial s}(z_0)-\frac{\partial f_2}{\partial p}(z_0)=e^{-i\theta}\frac{\partial f_2}{\partial s}(z_0)-e^{i\theta}\frac{\partial f_1}{\partial p}(z_0)$$ and $$e^{i\theta}\frac{\partial f_1}{\partial p}(z_0)+e^{-i\theta}\frac{\partial f_2}{\partial s}(z_0)=0.$$
\par Let $\Delta$ be the discriminant of the characteristic polynomial of $J_f(z_0)$, then
\begin{align*}
\Delta&=\left(\frac{\partial f_1}{\partial s}(z_0)+\frac{\partial f_2}{\partial p}(z_0)\right)^2-4\left(\frac{\partial f_1}{\partial s}(z_0)\frac{\partial f_2}{\partial p}(z_0)-\frac{\partial f_1}{\partial p}(z_0)\frac{\partial f_2}{\partial s}(z_0)\right)\\
&=\left(\frac{\partial f_1}{\partial s}(z_0)-\frac{\partial f_2}{\partial p}(z_0)\right)^2+4\frac{\partial f_1}{\partial p}(z_0)\frac{\partial f_2}{\partial s}(z_0)\\
&=\left(e^{-i\theta}\frac{\partial f_2}{\partial s}(z_0)-e^{i\theta}\frac{\partial f_1}{\partial p}(z_0)\right)^2+4\frac{\partial f_1}{\partial p}(z_0)\frac{\partial f_2}{\partial s}(z_0)\\
&=\left(e^{-i\theta}\frac{\partial f_2}{\partial s}(z_0)+e^{i\theta}\frac{\partial f_1}{\partial p}(z_0)\right)^2\\
&=0.
\end{align*}
Thus, $\lambda=\mu$ is a root of order 2 of the characteristic polynomial of $J_f(z_0)$.
\par Therefore, $\lambda$ and $\mu$ are the all eigenvalues of the linear transformation $J_f(z_0)$ on $\mathbb{C}^2$. So, from $|\mu|\leqslant1$ by $\mathrm{(\romannumeral4)}$,  we get $$|\det J_f(z_0)|\leqslant\lambda,\,\,\,|\mathrm{tr}\,J_f(z_0)|\leqslant\lambda+1.$$
The proof of $\mathrm{(\romannumeral1)}$ and $\mathrm{(\romannumeral5)}$ is complete.
\par Finally we show that the inequalities in $\mathrm{(\romannumeral4)}$ and $\mathrm{(\romannumeral5)}$ are sharp. Obviously the identity mapping on $\emph{\textbf{G}}_2$ is an example to make the inequalities in $\mathrm{(\romannumeral4)}$ and $\mathrm{(\romannumeral5)}$ as equalities. Beside the identity mapping, we can also consider the holomorphic mapping
$$f(s,p)=\big(e^{-i\theta}\frac{e^{i\theta}s-2p}{2-e^{-i\theta}s}+e^{-i\theta}p, e^{-2i\theta}\frac{e^{i\theta}s-2p}{2-e^{-i\theta}s}p\big):\emph{\textbf{G}}_2 \rightarrow \emph{\textbf{G}}_2.$$
Then we obtain $\mu=\frac{\partial f_1}{\partial s}(z_0)+e^{i\theta}\frac{\partial f_1}{\partial p}(z_0)=1$. That means the inequalities in $\mathrm{(\romannumeral4)}$ and $\mathrm{(\romannumeral5)}$ are sharp. The proof is complete.

\begin{thm} \label{Th:3.2}
Let $f:\emph{\textbf{G}}_2\rightarrow\emph{\textbf{G}}_2$ be a holomorphic mapping with $f(0)=0$ and set $f=(f_1, f_2)$. Let $z_0=(0,e^{i\theta})\in\partial\emph{\textbf{G}}_2$, where $\theta\in[0,2\pi)$. If $f$ is holomorphic at $z_0$ and $f(z_0)=z_0$, then the following statements hold:

\par $\mathrm{(\romannumeral1)}$   $\lambda:=\frac{\partial f_2}{\partial p}(z_0)$ and  $\mu:=\frac{\partial f_1}{\partial s}(z_0)$ are all eigenvalues of $J_f(z_0)$;
\par $\mathrm{(\romannumeral2)}$   $\lambda=\frac{\partial f_2}{\partial p}(z_0)\geqslant1;$
\par $\mathrm{(\romannumeral3)}$   $\frac{\partial f_1}{\partial p}(z_0)=0$. That is, $J_f(z_0)$ is a lower triangular square marix, and $(0,1)'$ is an eigenvector of $J_f(z_0)$ with respect to $\lambda$;
\par $\mathrm{(\romannumeral4)}$   $\mu=\frac{\partial f_1}{\partial s}(z_0)$ satisfies $|\mu|\leqslant \lambda$, and $|\frac{\partial f_2}{\partial s}(z_0)|\leqslant\lambda;$
\par $\mathrm{(\romannumeral5)}$   $|\det J_f(z_0)|\leqslant\lambda^2$,   $|\mathrm{tr}\, J_f(z_0)|\leqslant 2\lambda.$
\par Moreover, the inequalities in $\mathrm{(\romannumeral2)}$, $\mathrm{(\romannumeral4)}$ and $\mathrm{(\romannumeral5)}$ are sharp.
\end{thm}
\noindent
{\bf Proof}.
The proof is divided into three steps.\\
\emph{Step} 1.  Take $\phi_1(\xi)=e^{-i\theta}f_2(0,\xi e^{i\theta})$, $\xi\in D$. Then $\phi_1:D\rightarrow D$ is holomorphic in $D\cup\{1\}$ and such that $\phi_1(0)=0$, $\phi_1(1)=e^{-i\theta}f_2(0,e^{i\theta})=1$. Then by Theorem \ref{Th:1.2}, we have
$$\phi_1'(1)=\frac{\partial f_2}{\partial p}(z_0)\geqslant1.$$
\par Let $$\phi_2(\xi)=e^{-i\theta}\frac{2f_2(0,\xi e^{i\theta})-\omega f_1(0,\xi e^{i\theta})}{2-\overline{\omega}f_1(0,\xi e^{i\theta})},\,\,\xi\in D,$$
where $\omega\in\mathbb{T}$ is any fixed complex number. Then from the equivalence of (1) and (4) in Lemma \ref{Lm:2.1}, we can easily get that
 $\phi_2:D\rightarrow D$ is holomorphic in $D\cup\{1\}$ and such that $\phi_2(0)=0$, $\phi_2(1)=e^{i\theta}\frac{2f_2(0,e^{i\theta})-\omega f_1(0,e^{i\theta})}{2-\overline{\omega}f_1(0,e^{i\theta})}=1$.
\par By Theorem \ref{Th:1.2}, we have $$\phi_2'(1)=\frac{\partial f_2}{\partial p}(z_0)-\frac{1}{2}(\omega-\overline{\omega}e^{i\theta})\frac{\partial f_1}{\partial p}(z_0)\geqslant1.$$
That implies $$(\omega-\overline{\omega}e^{i\theta})\frac{\partial f_1}{\partial p}(z_0)\in \mathbb{R}$$
 and
$$2\frac{\partial f_2}{\partial p}(z_0)\geqslant 2+(\omega-\overline{\omega}e^{i\theta})\frac{\partial f_1}{\partial p}(z_0)$$ for any $\omega\in\mathbb{T}$.
\par Let $\omega=a+bi$, where $a,b\in\mathbb{R}$. If $\theta=0$, then
$$(\omega-\overline{\omega}e^{i\theta})\frac{\partial f_1}{\partial p}(z_0)=(\omega-\overline{\omega})\frac{\partial f_1}{\partial p}(z_0)=2bi\frac{\partial f_1}{\partial p}(z_0)\in\mathbb{R},\,\,\,\forall b:\,\,-1 \leqslant b \leqslant 1.$$
Let $\frac{\partial f_1}{\partial p}(z_0)=ci$ with $c\in\mathbb{R}$. Suppose that $c\neq0$, then we have
$$2\frac{\partial f_2}{\partial p}(z_0)\geqslant2+(\omega-\overline{\omega})\frac{\partial f_1}{\partial p}(z_0)=2-2bc,\,\,\,\forall b:\,\,-1 \leqslant b \leqslant 1.$$
Take $b_0\in[-1,1]$ such that $b_0c<0$, then$$\frac{\partial f_2}{\partial p}(z_0)\geqslant1-b_0c>1.$$
That is a contradiction if we take the identity mapping on the $\emph{\textbf{G}}_2$. Hence we have $c=0$. It implies $\frac{\partial f_1}{\partial p}(z_0)=0$.\\
\par If $\theta\neq0$, then take $\omega=1$ and $\omega=e^{i\theta}$, we get $$\frac{\partial f_2}{\partial p}(z_0)\geqslant1+\frac{1}{2}(1-e^{i\theta})\frac{\partial f_1}{\partial p}(z_0)$$
and
$$\frac{\partial f_2}{\partial p}(z_0)\geqslant1+\frac{1}{2}(e^{i\theta}-1)\frac{\partial f_1}{\partial p}(z_0)$$
respectively. Suppose $\frac{\partial f_1}{\partial p}(z_0)\neq0$, then $(1-e^{i\theta})\frac{\partial f_1}{\partial p}(z_0)\neq0$.
 Without loss of generality, we can assume $(1-e^{i\theta})\frac{\partial f_1}{\partial p}(z_0)>0$. Consequently
 $$\frac{\partial f_2}{\partial p}(z_0)\geqslant1+\frac{1}{2}(1-e^{i\theta})\frac{\partial f_1}{\partial p}(z_0)>1.$$
Similar to the case of $\theta=0$ we know that it is a contradiction. Thus we have $\frac{\partial f_1}{\partial p}(z_0)=0$. That is, $J_f(z_0)$ is a lower triangular square marix, $\lambda=\frac{\partial f_2}{\partial p}(z_0)\geqslant1$ is an eigenvalue of $J_f(z_0)$, $(0,1)'$ is an eigenvector of $J_f(z_0)$ with respect to $\lambda$.
\par The proof of $\mathrm{(\romannumeral1)}$, $\mathrm{(\romannumeral2)}$ and $\mathrm{(\romannumeral3)}$ is complete.
\\

\emph{Step} 2.  By Lemma \ref{Lm:2.2}, for any $p\in\mathbb{C}$ with $0<|p|<1$, we have
$$F_C((0,p),(\xi_1,0)')=\sup_{\omega\in\mathbb{T}}\frac{|\xi_1(1-\omega^2 p)|}{2(1-|p|^2)}=\frac{|\xi_1|}{2(1-|p|^2)}\sup_{\omega\in\mathbb{T}}|1-\omega^2p|=\frac{|\xi_1|}{2(1-|p|)}.$$
By Lemma \ref{Lm:2.4}, for any $0<t<1$, take any fixed $\xi_1\in\mathbb{C}$ with $\xi_1\neq0$, we have
\begin{equation}\label{eq3.2}
F_C\big(f(0,te^{i\theta}),J_f(0,te^{i\theta})(\xi_1,0)'\big)\leqslant F_C((0,te^{i\theta}),(\xi_1,0)').
\end{equation}
It follows that
\begin{align}
\nonumber & \;\;\;\;\;  \frac{|\xi_1|}{2(1-t)}\\
\nonumber & =F_C((0,te^{i\theta}),(\xi_1,0)')\\
\nonumber & \geqslant F_C\big(f(0,te^{i\theta}),J_f(0,te^{i\theta})(\xi_1,0)'\big)\\
\nonumber&=2\sup_{\omega\in\mathbb{T}}\frac{\bigg|J_{f_1}(0,te^{i\theta})\left( \begin{array}{ccc}
\xi_1 \\ 0  \\
 \end{array} \right)(1-\omega^2f_2(0,te^{i\theta}))-J_{f_2}(0,te^{i\theta})\left( \begin{array}{ccc}
\xi_1 \\ 0  \\
 \end{array} \right)(2-\omega f_1(0,te^{i\theta}))\omega\bigg|}{|2-\omega f_1(0,te^{i\theta})|^2-|2\omega f_2(0,te^{i\theta})-f_1(0,te^{i\theta})|^2}\\
 \nonumber&\geqslant2\frac{\bigg|J_{f_1}(0,te^{i\theta})\left( \begin{array}{ccc}
\xi_1 \\ 0  \\
 \end{array} \right)(1-e^{-i\theta}f_2(0,te^{i\theta}))-J_{f_2}(0,te^{i\theta})\left( \begin{array}{ccc}
\xi_1 \\ 0  \\
 \end{array} \right)(2-e^{-\frac{\theta}{2}i} f_1(0,te^{i\theta}))e^{-\frac{\theta}{2}i}\bigg|}{|2-e^{-\frac{\theta}{2}i} f_1(0,te^{i\theta})|^2-|2e^{-\frac{\theta}{2}i} f_2(0,te^{i\theta})-f_1(0,te^{i\theta})|^2}.
\end{align}
Notice that $$f_1(0,te^{i\theta})=o(t-1)$$
and $$f_2(0,te^{i\theta})=e^{i\theta}+e^{i\theta}\lambda(t-1)+o(t-1),$$
here $t$ is in a left neighborhood of $1$.
Thus
$$|2-e^{-\frac{\theta}{2}i}f_1(0,te^{i\theta})|^2-|2e^{-\frac{\theta}{2}i}f_2(0,te^{i\theta})-f_1(0,te^{i\theta})|^2=8\lambda(1-t)+o(t-1).$$
This means
\begin{equation*}
\begin{aligned}
\frac{|\xi_1|}{2}\frac{8\lambda(1-t)+o(t-1)}{1-t}\geqslant&2\bigg|J_{f_1}(0,te^{i\theta})\left( \begin{array}{ccc}
\xi_1 \\ 0  \\
 \end{array} \right)(1-e^{-i\theta}f_2(0,te^{i\theta})) \\ &  \;\; -J_{f_2}(0,te^{i\theta})\left( \begin{array}{ccc}
\xi_1 \\ 0  \\
 \end{array} \right)(2-e^{-\frac{\theta}{2}i} f_1(0,te^{i\theta}))e^{-\frac{\theta}{2}i}\bigg|.
\end{aligned}
\end{equation*}
As $t\rightarrow1^-$, then we obtain $|\frac{\partial f_2}{\partial s}(z_0)|\leqslant\lambda.$
\\

\emph{Step} 3.  By Lemma \ref{Lm:2.2}, for any $0<t<1$ and $(\xi_1,\xi_2)'\in\mathbb{C}^2$, we have
$$F_C((0,te^{i\theta}),(\xi_1,\xi_2)')=\frac{1}{2(1-t^2)}\sup_{\omega\in \mathbb{T}}\big|\xi_1(1-\omega^{2}te^{i\theta} )-2\omega \xi_2\big|.$$
\par One can easily prove that $\big|\xi_1(1-\omega^{2}te^{i\theta} )-2\omega \xi_2\big|$ is uniformly continuous on $(t,\omega)\in[0,1]\times\mathbb{T}$, it follows that
$$\lim_{t\rightarrow1^{-}}\sup_{\omega\in\mathbb{T}}\big|\xi_1(1-\omega^{2}te^{i\theta} )-2\omega \xi_2\big|=\sup_{\omega\in\mathbb{T}}\big|\xi_1(1-\omega^{2}e^{i\theta} )-2\omega \xi_2\big|,$$
so we have
$$\lim_{t\rightarrow1^{-}}(1-t)F_C((0,te^{i\theta}),(\xi_1,\xi_2)')=\frac{1}{4}\sup_{\omega\in\mathbb{T}}\big|\xi_1(1-\omega^{2}e^{i\theta} )-2\omega \xi_2\big|.$$
\par Similarly, we can get
\begin{align}
\nonumber & \;\;\;\;\; \lim_{t\rightarrow1^{-}}(1-t)F_C(f(0,te^{i\theta}),\emph{J}_f(0,te^{i\theta})(\xi_1,\xi_2)') \\
\nonumber & =\frac{1}{4\lambda}\sup_{\omega\in\mathbb{T}}\bigg|J_{f_1}(z_0)\left( \begin{array}{ccc}
\xi_1 \\ \xi_2  \\
 \end{array} \right)(1-\omega^{2}e^{i\theta} )-2J_{f_2}(z_0)\left( \begin{array}{ccc}
\xi_1 \\ \xi_2  \\
 \end{array} \right)\omega \bigg|.
\end{align}
\par Notice that
$$F_C((0,te^{i\theta}),(\xi_1,\xi_2)')\geqslant F_C(f(0,te^{i\theta}),\emph{J}_f(0,te^{i\theta})(\xi_1,\xi_2)'),$$
so we have
$$(1-t)F_C((0,te^{i\theta}),(\xi_1,\xi_2)')\geqslant (1-t)F_C(f(0,te^{i\theta}),\emph{J}_f(0,te^{i\theta})(\xi_1,\xi_2)').$$
As $t\rightarrow1^-$, we obtain
$$\lambda\sup_{\omega\in\mathbb{T}}\big|\xi_1(1-\omega^{2}e^{i\theta} )-2\omega \xi_2\big|\geqslant\sup_{\omega\in\mathbb{T}}\bigg|J_{f_1}(z_0)\left( \begin{array}{ccc}
\xi_1 \\ \xi_2  \\
 \end{array} \right)(1-\omega^{2}e^{i\theta} )-2J_{f_2}(z_0)\left( \begin{array}{ccc}
\xi_1 \\ \xi_2  \\
 \end{array} \right)\omega \bigg|.$$
\par Let $(\xi_1,\xi_2)'$ be an eigenvector of $J_f(z_0)$ with respect to the eigenvalue $\mu=\frac{\partial f_1}{\partial s}(z_0)$ , then it follows that

$$\lambda\sup_{\omega\in\mathbb{T}}\big|\xi_1(1-\omega^{2}e^{i\theta} )-2\omega \xi_2\big|\geqslant\big|\frac{\partial f_1}{\partial s}(z_0)\big|\sup_{\omega\in\mathbb{T}}\big|\xi_1(1-\omega^{2}e^{i\theta} )-2\omega \xi_2\big|.$$

\par We next prove that $$\sup_{\omega\in\mathbb{T}}\big|\xi_1(1-\omega^{2}e^{i\theta} )-2\omega \xi_2\big|\neq0.$$ In fact, assume $$\sup_{\omega\in\mathbb{T}}\big|\xi_1(1-\omega^{2}e^{i\theta} )-2\omega \xi_2\big|=0,$$
that is, for any $\omega\in\mathbb{T}$, we have $$\xi_1(1-\omega^{2}e^{i\theta} )-2\omega \xi_2=0.$$
Take $\omega =e^{-\frac{\theta}{2}i}$, then we get $\xi_2=0$. We can also take $\omega=1$ and obtain $\xi_1=0$. So we get $\xi_1=\xi_2=0$, which is in contradiction to the fact that $(\xi_1,\xi_2)'$ is an eigenvector of $J_f(z_0)$ with respect to the eigenvalue $\mu=\frac{\partial f_1}{\partial s}(z_0)$. This means $$\sup_{\omega\in\mathbb{T}}\big|\xi_1(1-\omega^{2}e^{i\theta} )-2\omega \xi_2\big|\neq0.$$ It follows that $$|\mu|=\big|\frac{\partial f_1}{\partial s}(z_0)\big|\leqslant\lambda.$$
So we can easily get $$|\det J_f(z_0)|\leqslant\lambda^2,\,\,|\mathrm{tr}\, J_f(z_0)|\leqslant 2\lambda.$$
The proof of $\mathrm{(\romannumeral4)}$ and $\mathrm{(\romannumeral5)}$ is complete.
\par Finally we show that the inequalities in $\mathrm{(\romannumeral2)}$, $\mathrm{(\romannumeral4)}$ and $\mathrm{(\romannumeral5)}$ are sharp. Obviously the identity mapping on symmetrized bidisc $\emph{\textbf{G}}_2$ shows that these inequalities are sharp. Beside the identity mapping on $\emph{\textbf{G}}_2$, there are so many other examples. Considering the holomorphic mapping
$$f(s,p)=\big(0,-ie^{\frac{\theta}{2}i}\frac{2ie^{-\frac{\theta}{2}i}p-s}{2-ie^{-\frac{\theta}{2}i}s}\big):{\textbf{G}}_2 \rightarrow {\textbf{G}}_2.$$
Then we have $\lambda=\frac{\partial f_2}{\partial p}(z_0)=1$ and $\frac{\partial f_2}{\partial s}(z_0)=ie^{\frac{\theta}{2}i}$. That means the inequality in $\mathrm{(\romannumeral2)}$ and the inequality in $\mathrm{(\romannumeral4)}$ are sharp. In order to verify that the first inequality in $\mathrm{(\romannumeral4)}$ and the inequalities in $\mathrm{(\romannumeral5)}$ are sharp, we consider the holomorphic mapping
$$f(s,p)=\big(\frac{s-2\omega_1p}{2-\omega_1s}+\frac{s-2\omega_2p}{2-\omega_2s},\frac{s-2\omega_1p}{2-\omega_1s}\frac{s-2\omega_2p}{2-\omega_2s}\big):\emph{\textbf{G}}_2\rightarrow \emph{\textbf{G}}_2,$$
where $\omega_1=-ie^{-\frac{\theta}{2}i}$, $\omega_2=ie^{-\frac{\theta}{2}i}$. Then a simple calculation shows that $\mu=\frac{\partial f_1}{\partial s}(z_0)=2$ and $\lambda=\frac{\partial f_2}{\partial p}(z_0)=2$. That means the first inequality in $\mathrm{(\romannumeral4)}$ and the inequalities in $\mathrm{(\romannumeral5)}$ are sharp. The proof is complete.

\vskip 6pt
\noindent {\bf Remark.} The condition $f(0)=0$ in Theorem \ref{Th:3.2} can not be removed. Consider
$$H(s,p)=\bigg(\frac{(1+r^2)s+2rie^{-\frac{\theta}{2}i}p-2rie^{\frac{\theta}{2}i}}{1+rie^{-\frac{\theta}{2}i}s-r^2e^{-i\theta}p},\frac{p-rie^{\frac{\theta}{2}i}s-r^2e^{i\theta}}{1+rie^{-\frac{\theta}{2}i}s-r^2e^{-i\theta}p}\bigg),$$where $0<r<1$. Then $H(s,p):\emph{\textbf{G}}_2\rightarrow \emph{\textbf{G}}_2$ is a holomorphic mapping and $H(0,e^{i\theta})=(0,e^{i\theta})$, $H(0)\neq0$, but $\frac{\partial H}{\partial p}(z_0)=\frac{2ri}{1-r^2}e^{-\frac{\theta}{2}i}\neq0$. That meas $J_H(z_0)$ is not a lower triangular square marix.

\begin{thm} \label{Th:3.3}
Let $f:\emph{\textbf{G}}_2\rightarrow\emph{\textbf{G}}_2$ be a holomorphic mapping with $f(0)=0$ and set $f=(f_1, f_2)$. Let $z_0=(2\alpha,\alpha^2)\in\partial\emph{\textbf{G}}_2$ with $\alpha\in\mathbb{C}$ and $|\alpha|=1$. If $f$ is holomorphic at $z_0$ and $f(z_0)=z_0$, then for the eigenvalues $\lambda$, $\mu$ of $J_f(z_0)$, the following statements hold:
\par $\mathrm{(\romannumeral1)}$   $\lambda=\frac{\partial f_1}{\partial s}(z_0)+\alpha\frac{\partial f_1}{\partial p}(z_0)=
\overline{\alpha}\frac{\partial f_2}{\partial s}(z_0)+\frac{\partial f_2}{\partial p}(z_0)\geqslant1$;
\par $\mathrm{(\romannumeral2)}$   $(1,\alpha)'$ is an eigenvector of $J_f(z_0)$ with respect to $\lambda$. That is
$$J_f(z_0)
\left( \begin{array}{ccc}
1 \\ \alpha  \\
 \end{array} \right)=\lambda\left( \begin{array}{ccc}
1 \\ \alpha  \\
 \end{array} \right);$$
\par $\mathrm{(\romannumeral3)}$  Let
$$A=2\frac{\partial^2 f_1}{\partial s^2}(z_0)+4\alpha\frac{\partial^2 f_1}{\partial s \partial p}(z_0)+\frac{\partial f_1}{\partial p}(z_0)+2\alpha^2\frac{\partial^2 f_1}{\partial p^2}(z_0),$$
$$B=2\frac{\partial^2 f_2}{\partial s^2}(z_0)+4\alpha\frac{\partial^2 f_2}{\partial s \partial p}(z_0)+\frac{\partial f_2}{\partial p}(z_0)+2\alpha^2\frac{\partial^2 f_2}{\partial p^2}(z_0).$$
Then $B-A\alpha\in\mathbb{R}$ and $|\mu|\leqslant B-A\alpha$;
\par $\mathrm{(\romannumeral4)}$  $|\det J_f(z_0)|\leqslant(B-A\alpha)\lambda,\,\,\,|\mathrm{tr}\, J_f(z_0)|\leqslant \lambda+B-A\alpha.$
\par Moreover, the inequalities in $\mathrm{(\romannumeral1)}$, $\mathrm{(\romannumeral3)}$ and $\mathrm{(\romannumeral4)}$ are sharp.
\end{thm}
\noindent
{\bf Proof}.
The proof is divided into three steps.\\
\emph{Step} 1.   Take $\psi_1(\xi)=\frac{1}{2}\overline{\alpha}f_1(2\xi\alpha,\xi^2\alpha^2)$, $\xi\in D$. Then $\psi_1:D\rightarrow D$ is holomorphic in $D\cup\{1\}$ and such that $\psi_1(0)=0$, $\psi_1(1)=\frac{1}{2}\overline{\alpha}f_1(2\alpha,\alpha^2)=1$. Then by Theorem \ref{Th:1.2}, we have
$$\lambda_1:=\psi_1'(1)=\frac{\partial f_1}{\partial s}(z_0)+\alpha\frac{\partial f_1}{\partial p}(z_0)\geqslant1.$$
\par Take $\psi_2(\xi)=\overline{\alpha}^2f_2(2\xi\alpha,\xi^2\alpha^2)$, $\xi\in D$. Then $\psi_2:D\rightarrow D$ is holomorphic in $D\cup\{1\}$ and such that $\psi_2(0)=0$, $\psi_2(1)=\overline{\alpha}^2f_2(2\alpha,\alpha^2)=1$. Then by Theorem \ref{Th:1.2} again, we have
$$\psi_2'(1)=2\Big(\overline{\alpha}\frac{\partial f_2}{\partial s}(z_0)+\frac{\partial f_2}{\partial p}(z_0)\Big)\geqslant1.$$
Let $\lambda_2=\overline{\alpha}\frac{\partial f_2}{\partial s}(z_0)+\frac{\partial f_2}{\partial p}(z_0)$. Now we want to prove $\lambda_1=\lambda_2$. The proof is divided into three cases.
\par Case 1. $\alpha=1$.  For any fixed $\theta:0<\theta<\pi$, take
$$\psi_3(\xi)=e^{i\theta}\frac{e^{-i\theta} f_1(2\xi,\xi^2)-2f_2(2\xi,\xi^2)}{2-e^{i\theta}f_1(2\xi,\xi^2)},\,\,\,\xi\in D.$$
 From the equivalence of (1) and (4) in Lemma \ref{Lm:2.1}, we can easily get that $\psi_3:D\rightarrow D$ is holomorphic in $D\cup\{1\}$ and such that $\psi_3(0)=0$, $\psi_3(1)=1$. Notice that
$$\lambda_1=\frac{\partial f_1}{\partial s}(z_0)+\frac{\partial f_1}{\partial p}(z_0)\in\mathbb{R},\,\,\lambda_2=\frac{\partial f_2}{\partial s}(z_0)+\frac{\partial f_2}{\partial p}(z_0)\in\mathbb{R}.$$
Then by Theorem \ref{Th:1.2}, we have
\begin{align*}
\psi_3'(1)&=\frac{(1+e^{i\theta})\lambda_1-2e^{i\theta}\lambda_2}{1-e^{i\theta}}\\
&=\frac{(1+e^{i\theta})(1-e^{-i\theta})\lambda_1-2e^{i\theta}(1-e^{-i\theta})\lambda_2}{|1-e^{i\theta}|^2}\\
&=\frac{(e^{i\theta}-e^{-i\theta})\lambda_1-2(e^{i\theta}-1)\lambda_2}{|1-e^{i\theta}|^2}\geqslant1.
\end{align*}
It follows that $$(e^{i\theta}-e^{-i\theta})\lambda_1-2(e^{i\theta}-1)\lambda_2=2(1-\cos\theta)\lambda_2+2\sin\theta (\lambda_1-\lambda_2)i\in\mathbb{R},$$
which implies $\lambda_1=\lambda_2.$
\par Case 2.  $\alpha=-1$.  For any fixed $\theta:0<\theta<\pi$, take
$$\psi_3(\xi)=-e^{i\theta}\frac{e^{-i\theta}f_1(-2\xi,\xi^2)-2f_2(-2\xi,\xi^2)}{2-e^{i\theta}f_1(-2\xi,\xi^2)},\,\,\,\xi\in D.$$
Then $\psi_3:D\rightarrow D$ is is holomorphic in $D\cup\{1\}$ and such that $\psi_3(0)=0$, $\psi_3(1)=1$. Notice that
$$\lambda_1=\frac{\partial f_1}{\partial s}(z_0)-\frac{\partial f_1}{\partial p}(z_0)>0,\,\,\lambda_2=-\frac{\partial f_2}{\partial s}(z_0)+\frac{\partial f_2}{\partial p}(z_0)>0.$$ Then we have
\begin{align*}
\psi_3'(1)&=\frac{(1-e^{i\theta})\lambda_1+2e^{i\theta}\lambda_2}{1+e^{i\theta}}\\
&=\frac{(1-e^{i\theta})(1+e^{-i\theta})\lambda_1+2e^{i\theta}(1+e^{-i\theta})\lambda_2}{|1+e^{i\theta}|^2}\\
&=\frac{(e^{-i\theta}-e^{i\theta})\lambda_1+2(e^{i\theta}+1)\lambda_2}{|1+e^{i\theta}|^2}\geqslant1.
\end{align*}
It follows that $$(e^{-i\theta}-e^{i\theta})\lambda_1+2(1+e^{i\theta})\lambda_2=2(1+\cos\theta)\lambda_2+2\sin\theta(\lambda_2-\lambda_1)i\in\mathbb{R},$$
which implies $\lambda_1=\lambda_2.$
\par Case 3.  $\alpha\neq-1\,\, \mathrm{and}\,\, 1$.   Take
$$\psi_3(\xi)=\overline{\alpha}\frac{f_1(2\xi\alpha,\xi^2\alpha^2)-2f_2(2\xi\alpha,\xi^2\alpha^2)}{2-f_1(2\xi\alpha,\xi^2\alpha^2)},\,\,\,\xi\in D.$$
Then we obtain
\begin{align*}
\psi_3'(1)&=\frac{(1+\alpha)\lambda_1-2\alpha\lambda_2}{1-\alpha}\\
&=\frac{(1+\alpha)(1-\overline{\alpha})\lambda_1-2\alpha(1-\overline{\alpha})\lambda_2}{|1-\alpha|^2}\\
&=\frac{(\alpha-\overline{\alpha})\lambda_1-2(\alpha-1)\lambda_2}{|1-\alpha|^2}\geqslant1.
\end{align*}
This implies $(\alpha-\overline{\alpha})\lambda_1-2(\alpha-1)\lambda_2\in\mathbb{R}$. Let $\alpha=a+bi$. Since $\alpha\neq-1\,\, \mathrm{and}\,\, 1$, we have $b\neq0$. Notice that
$$(\alpha-\overline{\alpha})\lambda_1-2(\alpha-1)\lambda_2=2(1-a)\lambda_2+2b(\lambda_1-\lambda_2)i\in\mathbb{R},$$
it follows that $\lambda_1=\lambda_2$.
\par Let $\lambda:=\lambda_1=\lambda_2\geqslant1$, then
$$\lambda=\frac{\partial f_1}{\partial s}(z_0)+\alpha\frac{\partial f_1}{\partial p}(z_0)=\overline{\alpha}\frac{\partial f_2}{\partial s}(z_0)+\frac{\partial f_2}{\partial p}(z_0)\geqslant1.$$
It follows that
$$J_f(z_0)\left( \begin{array}{ccc}
1 \\ \alpha  \\
 \end{array} \right)=\left( \begin{array}{ccc}
\frac{\partial f_1}{\partial s}(z_0)+\alpha\frac{\partial f_1}{\partial p}(z_0) \\ \frac{\partial f_2}{\partial s}(z_0)+\alpha\frac{\partial f_2}{\partial p}(z_0)  \\
 \end{array} \right)=
\left( \begin{array}{ccc}
\lambda \\ \alpha\lambda  \\
 \end{array} \right)=\lambda
\left( \begin{array}{ccc}
1 \\ \alpha  \\
 \end{array} \right).$$
That means $\lambda$ is an eigenvalue of $J_f(z_0)$ and $(1,\alpha)'$ is an eigenvector of $J_f(z_0)$ with respect to $\lambda$.
\par The proof of $\mathrm{(\romannumeral1)}$ and $\mathrm{(\romannumeral2)}$ is complete.
\\ \\

\emph{Step} 2.  For $0<t<1$, consider the Taylor expansion of $f_1(2t\alpha,t^2\alpha^2)$ and $f_2(2t\alpha,t^2\alpha^2)$ at $t=1$, we have
\begin{equation*}
\begin{aligned}
f_1(2t\alpha,t^2\alpha^2)=&2\alpha+2\alpha\lambda(t-1)\\
&+\alpha^2\Big(2\frac{\partial^2 f_1}{\partial s^2}(z_0)+4\alpha\frac{\partial^2 f_1}{\partial s \partial p}(z_0)+\frac{\partial f_1}{\partial p}(z_0)+2\alpha^2\frac{\partial^2 f_1}{\partial p^2}(z_0)\Big)(t-1)^2+o((t-1)^2),
\end{aligned}
\end{equation*}
\begin{equation*}
\begin{aligned}
f_2(2t\alpha,t^2\alpha^2)=&\alpha^2+2\alpha^2\lambda(t-1)\\
&+\alpha^2\Big(2\frac{\partial^2 f_2}{\partial s^2}(z_0)+4\alpha\frac{\partial^2 f_2}{\partial s \partial p}(z_0)+\frac{\partial f_2}{\partial p}(z_0)+2\alpha^2\frac{\partial^2 f_2}{\partial p^2}(z_0)\Big)(t-1)^2+o((t-1)^2),
\end{aligned}
\end{equation*}
here $t$ is in a left neighborhood of $1$.
\par If
\begin{equation*}
A=2\frac{\partial^2 f_1}{\partial s^2}(z_0)+4\alpha\frac{\partial^2 f_1}{\partial s \partial p}(z_0)+\frac{\partial f_1}{\partial p}(z_0)+2\alpha^2\frac{\partial^2 f_1}{\partial p^2}(z_0),
\end{equation*}
\begin{equation*}
B=2\frac{\partial^2 f_2}{\partial s^2}(z_0)+4\alpha\frac{\partial^2 f_2}{\partial s \partial p}(z_0)+\frac{\partial f_2}{\partial p}(z_0)+2\alpha^2\frac{\partial^2 f_2}{\partial p^2}(z_0),
\end{equation*}
then we have
\begin{equation}\label{eq3.4}
f_1(2t\alpha,t^2\alpha^2)=2\alpha+2\alpha\lambda(t-1)+\alpha^2 A(t-1)^2+o((t-1)^2),
\end{equation}
\begin{equation}\label{eq3.5}
f_2(2t\alpha,t^2\alpha^2)=\alpha^2+2\alpha^2\lambda(t-1)+\alpha^2 B(t-1)^2+o((t-1)^2).
\end{equation}
It follows that
\begin{equation}\label{eq3.6}
|2-\overline{\alpha}f_1(2t\alpha,t^2\alpha^2)|^2-|2\overline{\alpha}f_2(2t\alpha,t^2\alpha^2)-f_1(2t\alpha,t^2\alpha^2)|^2=
8\lambda\Re (B-A\alpha)\,(1-t)^3+o((t-1)^3).
\end{equation}
\par From \eqref{eq3.4} and \eqref{eq3.5}, we also have
\begin{equation}\label{eq3.7}
1-\overline{\alpha}^2f_2(2t\alpha,t^2\alpha^2)=2\lambda(1-t)-B(1-t)^2+o((1-t)^2),
\end{equation}
\begin{equation}\label{eq3.8}
(2-\overline{\alpha}f_1(2t\alpha,t^2\alpha^2))\overline{\alpha}=2\overline{\alpha}\lambda(1-t)-A(1-t)^2+o((1-t)^2).
\end{equation}
By Lemma \ref{Lm:2.5}, for any $0<t<1$, we have
\begin{equation}\label{eq3.9}
F_C((2t\alpha,t^2\alpha^2),(\xi_1,\xi_2)')=\frac{|(1+t^2)\xi_1-2t\overline{\alpha}\xi_2|+2|t\alpha\xi_1-\xi_2|}{2(1-t^2)^2}.
\end{equation}
This, together with \eqref{eq3.6}, \eqref{eq3.7} and \eqref{eq3.8}, implies that for any $0<t<1$,
\begin{align*}
&\,\,\,\,\,\,\,\frac{3-t}{2(1-t)(1+t)^2}\\&=F_C((2t\alpha,t^2\alpha^2),(1,\alpha)')\\
&\geqslant F_C\Big[f(2t\alpha,t^2\alpha^2),J_f(2t\alpha,t^2\alpha^2)\left( \begin{array}{ccc}
1 \\ \alpha  \\
 \end{array} \right)\Big]\\
&\geqslant2\,\,\frac{\bigg|J_{f_1}(2t\alpha,t^2\alpha^2)\left( \begin{array}{ccc}
1 \\ \alpha  \\
 \end{array} \right)(1-\overline{\alpha}^2f_2(2t\alpha,t^2\alpha^2))-J_{f_2}(2t\alpha,t^2\alpha^2)\left( \begin{array}{ccc}
1 \\ \alpha  \\
 \end{array} \right)(2-\overline{\alpha} f_1(2t\alpha,t^2\alpha^2))\overline{\alpha}\bigg|}{|2-\overline{\alpha} f_1(2t\alpha,t^2\alpha^2)|^2-|2\overline{\alpha} f_2(2t\alpha,t^2\alpha^2)-f_1(2t\alpha,t^2\alpha^2)|^2}\\
&=2\,\,\frac{|(B-A\alpha)\lambda(1-t)^2+o((1-t)^2)|}{8\lambda\Re(B-A\alpha)\,(1-t)^3+o((1-t)^3)}.
\end{align*}
That is
$$\frac{3-t}{2(1+t)^2}\geqslant\frac{|(B-A\alpha)\lambda(1-t)^3+o((1-t)^3)|}{4\lambda\Re(B-A\alpha)\,(1-t)^3+o((1-t)^3)}.$$
As $t\rightarrow1^-$, we obtain
$$|B-A\alpha|\leqslant\Re(B-A\alpha).$$
That means $B-A\alpha\in\mathbb{R}$. Moreover, we have $B-A\alpha\geqslant0$. The proof of the first part of $\mathrm{(\romannumeral4)}$ is complete.
\\ \\

\emph{Step} 3. From \eqref{eq3.9}, together with the result in step 2, we obtain
\begin{align*}
&\,\,\,\,\,\,\,\frac{|(1+t^2)\xi_1-2t\overline{\alpha}\xi_2|+2|t\alpha\xi_1-\xi_2|}{2(1-t^2)^2}\\
&=F_C((2t\alpha,t^2\alpha^2),(\xi_1,\xi_2)')\\
&\geqslant F_C(f(2t\alpha,t^2\alpha^2),J_f(2t\alpha,t^2\alpha^2)(\xi_1,\xi_2)')\\
&\geqslant2\frac{\bigg|J_{f_1}(2t\alpha,t^2\alpha^2)\left( \begin{array}{ccc}
\xi_1 \\ \xi_2  \\
 \end{array} \right)(1-\overline{\alpha}^2f_2(2t\alpha,t^2\alpha^2))-J_{f_2}(2t\alpha,t^2\alpha^2)\left( \begin{array}{ccc}
\xi_1 \\ \xi_2  \\
 \end{array} \right)(2-\overline{\alpha} f_1(2t\alpha,t^2\alpha^2))\overline{\alpha}\bigg|}{8\lambda(B-A\alpha)(1-t)^3+o((1-t)^3)}.
\end{align*}
Take $\xi_1=-\frac{\partial f_1}{\partial p}(z_0)$ and $\xi_2=\frac{\partial f_1}{\partial s}(z_0)$. As $t\rightarrow1^-$, then we obtain
$$|\det J_f(z_0)|=\big|\frac{\partial f_1}{\partial s}(z_0)\frac{\partial f_2}{\partial p}(z_0)-\frac{\partial f_2}{\partial s}(z_0)\frac{\partial f_1}{\partial p}(z_0)\big|\leqslant(B-A\alpha)\lambda.$$
Since $\det J_f(z_0)=\mu\lambda$, it follows that
$$|\mu|\leqslant B-A\alpha.$$
Then we can conclude that $$|\mathrm{tr}\, J_f(z_0)|\leqslant \lambda+B-A\alpha.$$
The proof of $\mathrm{(\romannumeral3)}$ and $\mathrm{(\romannumeral4)}$ is complete.
\par Moreover, considering the identity mapping on symmetrized bidisc $\emph{\textbf{G}}_2$, we can find that the inequalities in $\mathrm{(\romannumeral1)}$, $\mathrm{(\romannumeral3)}$ and $\mathrm{(\romannumeral4)}$ are sharp. The proof is complete.
\vskip 6pt
\par From our main results Theorem \ref{Th:3.1}, Theorem \ref{Th:3.2} and Theorem \ref{Th:3.3}, we can obtain a boundary Schwarz lemma for holomorphic function $h(s,p):\emph{\textbf{G}}_2\rightarrow D$ and holomorphic mapping $\varphi(z):D\rightarrow\emph{\textbf{G}}_2$  at boundary points $(e^{i\theta},0)$, $(0,e^{i\theta})$ and $(2\alpha,\alpha^2)$ ($\theta\in[0,2\pi)$, $\alpha\in\mathbb{C}$ and $|\alpha|=1$) of $\emph{\textbf{G}}_2$, which we can summarize as following corollaries.
\begin{cor}\label{Co:3.1}
Let $h(s,p):\emph{\textbf{G}}_2\rightarrow D$ be a holomorphic function and let $z_0=(e^{i\theta},0)\in\partial\emph{\textbf{G}}_2$. If $h$ is holomorphic at $z_0$ and $h(z_0)\in\partial D$, suppose $h(z_0)=\omega e^{i\theta}$, where $\omega\in\mathbb{T}$. Then the following statements hold:
\par $\mathrm{(\romannumeral1)}$   $\frac{\partial h}{\partial s}(z_0)+e^{i\theta}\frac{\partial h}{\partial p}(z_0)=0$;
\par $\mathrm{(\romannumeral2)}$   $\overline{\omega}\frac{\partial h}{\partial s}(z_0)$ is a real number and
$$\overline{\omega}\frac{\partial h}{\partial s}(z_0)\geqslant\frac{1}{2}\frac{|1-\overline{\varphi(0)}|^2}{1-|\varphi(0)|^2}>0,$$
where $\varphi(0)=\frac{\overline{\omega}e^{-i\theta}h(0,0)}{2-\overline{\omega}e^{-i\theta}h(0,0)}$.
\end{cor}
\noindent
{\bf Proof}.
Define $f(s,p)=(\overline{\omega}h(s,p),0)$, then $f:\emph{\textbf{G}}_2\rightarrow\emph{\textbf{G}}_2$ is a holomorphic mapping and satisfies the conditions of Theorem \ref{Th:3.1}. Apply Theorem \ref{Th:3.1} to obtain the statements in the corollary.
\begin{cor}\label{Co:3.2}
Let $h(s,p):\emph{\textbf{G}}_2\rightarrow D$ be a holomorphic function with $h(0,0)=0$. Let $z_0=(0,e^{i\theta})\in\partial\emph{\textbf{G}}_2$. If $h$ is holomorphic at $z_0$ and $h(z_0)\in\partial D$, suppose $h(z_0)=\omega e^{i\theta}$, where $\omega\in\mathbb{T}$. Then the following statements hold:
\par $\mathrm{(\romannumeral1)}$   $\overline{\omega}\frac{\partial h}{\partial p}(z_0)$ is a real number and $\overline{\omega}\frac{\partial h}{\partial p}(z_0)\geqslant1$;
\par $\mathrm{(\romannumeral2)}$   $\big|\frac{\partial h}{\partial s}(z_0)\big|\leqslant\overline{\omega}\frac{\partial h}{\partial p}(z_0)$.
\par Moreover, the inequalities in $\mathrm{(\romannumeral1)}$ and $\mathrm{(\romannumeral2)}$ are sharp.
\end{cor}
\noindent
{\bf Proof}.
Define $f(s,p)=(0,\overline{\omega}h(s,p))$, then $f:\emph{\textbf{G}}_2\rightarrow\emph{\textbf{G}}_2$ is a holomorphic mapping and satisfies the conditions of Theorem \ref{Th:3.2}. From Theorem \ref{Th:3.2} we can obtain the statements in the corollary. Moreover, considering the mapping
$$\Phi(s,p)=\frac{2ie^{-\frac{\theta}{2}i}p-s}{2-ie^{-\frac{\theta}{2}i}s}:\emph{\textbf{G}}_2\rightarrow D,$$
then we can easily prove that the inequalities in corollary are sharp.
\begin{cor}\label{Co:3.3}
Let $h(s,p):\emph{\textbf{G}}_2\rightarrow D$ be a holomorphic function with $h(0,0)=0$. Let $z_0=(2\alpha,\alpha^2)\in\partial\emph{\textbf{G}}_2$, where $|\alpha|=1$. If $h$ is holomorphic at $z_0$ and $h(z_0)\in\partial D$, suppose $h(z_0)=\omega\alpha$, where $\omega\in\mathbb{T}$. Then
$$\overline{\omega}\big(\frac{\partial h}{\partial s}(z_0)+\alpha\frac{\partial h}{\partial p}(z_0)\big)\geqslant\frac{1}{2}.$$
\par Moreover, the inequality above is sharp.
\end{cor}
\noindent
{\bf Proof}.
Define $f(s,p)=(2\overline{\omega}h(s,p),\overline{\omega}^2h(s,p)^2)$, then $f:\emph{\textbf{G}}_2\rightarrow\emph{\textbf{G}}_2$ is a holomorphic mapping and satisfies the conditions of Theorem \ref{Th:3.3}. It follows from Theorem \ref{Th:3.3} that the Corollary is proved. Moreover, considering the mapping
$$\Phi(s,p)=\frac{s-2\beta p}{2-\beta s}:\emph{\textbf{G}}_2\rightarrow D,$$
where $\beta\in\mathbb{T}$ and $\beta\neq\overline{\alpha}$, then the inequality turns out to be sharp.

\vskip 6pt
\noindent {\bf Remark.} In fact, the above three corollaries can also be conducted from the Carath\'eodory metric of $\emph{\textbf{G}}_2$ and the boundary Schwarz lemma of the unit disc $D$, but the calculation process is very complicated. Here we give a simpler proof.
\begin{cor}\label{Co:3.4}
Let $\varphi(z)=(\varphi_1(z),\varphi_2(z)):D\rightarrow\emph{\textbf{G}}_2$ be a holomorphic mapping. Let $\lambda\in\partial D$ and $z_0=(e^{i\theta},0)\in\partial\emph{\textbf{G}}_2$. If $\varphi$ is holomorphic at $\lambda\in\partial D$ and $\varphi(\lambda)=z_0$, then
$$\lambda e^{-i\theta}(\varphi'_1(\lambda)-e^{-i\theta}\varphi'_2(\lambda))\geqslant\frac{1}{4}\frac{|1-\overline{g(0)}|^2}{1-|g(0)|^2}>0,$$
where $g(0)=e^{-2i\theta}\frac{e^{i\theta}\varphi_1(0)-2\varphi_2(0)}{2-e^{-i\theta}\varphi_1(0)}$;
\end{cor}
\noindent
{\bf Proof}.
Define
$$f(s,p)=\bigg(\varphi_1\big(\lambda e^{-i\theta}\frac{s-2e^{-i\theta}p}{2-e^{-i\theta}s}\big),\varphi_2\big(\lambda e^{-i\theta}\frac{s-2e^{-i\theta}p}{2-e^{-i\theta}s}\big)\bigg),$$
then $f:\emph{\textbf{G}}_2\rightarrow\emph{\textbf{G}}_2$ is a holomorphic mapping and satisfies the conditions of Theorem \ref{Th:3.1}. Apply Theorem \ref{Th:3.1} to obtain the corollary.
\begin{cor}\label{Co:3.5}
Let $\varphi(z)=(\varphi_1(z),\varphi_2(z)):D\rightarrow\emph{\textbf{G}}_2$ be a holomorphic mapping with $\varphi(0)=(0,0)$. Let $\lambda\in\partial D$ and $z_0=(0,e^{i\theta})\in\partial\emph{\textbf{G}}_2$. If $\varphi$ is holomorphic at $\lambda\in\partial D$ and $\varphi(\lambda)=z_0$, then the following statements hold:
\par $\mathrm{(\romannumeral1)}$ $\varphi_1'(\lambda)=0$;
\par $\mathrm{(\romannumeral2)}$ $e^{-i\theta}\lambda\varphi_2'(\lambda)$ is a real number and $e^{-i\theta}\lambda\varphi_2'(\lambda)\geqslant1$.
\par Moreover, the inequality in $\mathrm{(\romannumeral2)}$ is sharp.
\end{cor}
\noindent
{\bf Proof}.
Define
$$f(s,p)=\bigg(\varphi_1\big(\frac{2\omega p-s}{2-\omega s}\big),\varphi_2\big(\frac{2\omega p-s}{2-\omega s}\big)\bigg),$$
where $\omega=e^{-i\theta}\lambda$. Then $f:\emph{\textbf{G}}_2\rightarrow\emph{\textbf{G}}_2$ is a holomorphic mapping and satisfies the conditions of Theorem \ref{Th:3.2}. From Theorem \ref{Th:3.2} we can obtain the statements in the corollary. Moreover, for any fixed $\lambda\in\partial D$, considering the holomorphic mapping
$$\varphi(z)=(0,\overline{\lambda}e^{i\theta}z):D\rightarrow\emph{\textbf{G}}_2,$$
then we can easily prove that the inequality in $\mathrm{(\romannumeral2)}$ is sharp.
\begin{cor}\label{Co:3.6}
Let $\varphi(z)=(\varphi_1(z),\varphi_2(z)):D\rightarrow\emph{\textbf{G}}_2$ be a holomorphic mapping with $\varphi(0)=(0,0)$. Let $\lambda\in\partial D$ and $z_0=(2\alpha,\alpha^2)\in\partial\emph{\textbf{G}}_2$, where $|\alpha|=1$. If $\varphi$ is holomorphic at $\lambda\in\partial D$ and $\varphi(\lambda)=z_0$, then the following statements hold:
\par $\mathrm{(\romannumeral1)}$ $\varphi_1'(\lambda)-\overline{\alpha}\varphi_2'(\lambda)=0$;
\par $\mathrm{(\romannumeral2)}$ $\lambda\overline{\alpha}\varphi_1'(\lambda)=\lambda\overline{\alpha}^2\varphi_2'(\lambda)$ is a real number and $\lambda\overline{\alpha}\varphi_1'(\lambda)=\lambda\overline{\alpha}^2\varphi_2'(\lambda)\geqslant2$.
\par Moreover, the inequality in $\mathrm{(\romannumeral2)}$ is sharp.
\end{cor}
\noindent
{\bf Proof}.
Define
$$f(s,p)=\bigg(\varphi_1\big(\lambda\overline{\alpha}\frac{s-2\omega p}{2-\omega s}\big),\varphi_2\big(\lambda\overline{\alpha}\frac{s-2\omega p}{2-\omega s}\big)\bigg),$$
where $\omega$ is any fixed number with $|\omega|<1$. Then $f:\emph{\textbf{G}}_2\rightarrow\emph{\textbf{G}}_2$ is a holomorphic mapping and satisfies the conditions of Theorem \ref{Th:3.3}. It follows from Theorem \ref{Th:3.3} that the corollary is proved. Moreover, for any fixed $\lambda\in\partial D$, considering the holomorphic mapping
$$\varphi(z)=(2\overline{\lambda}\alpha z,\overline{\lambda}^2\alpha^2z^2):D\rightarrow\emph{\textbf{G}}_2,$$
then the inequality turns out to be sharp.

\vskip 10pt

\noindent\textbf{Acknowledgments}\quad The authors thank the referees for many useful comments.
 The project is supported by the National Natural Science Foundation of China (No. 11671306).

\end{document}